\documentclass[11pt]{article}

\usepackage{amsmath}
\usepackage{graphicx}
\usepackage{geometry}
\usepackage{amsthm}
\usepackage{bm}
\usepackage{hyperref}
\usepackage{amsfonts}
\usepackage{caption}
\usepackage{subcaption}
\usepackage{authblk}

\usepackage{appendix}

\hypersetup{
    colorlinks=true,
    linkcolor=magenta,
    filecolor=magenta,
    urlcolor=magenta,
}

\newtheorem{thm}{Theorem}

\title{Arbitrarily high order implicit ODE integration by correcting a neural network approximation with Newton's method}
\date{\today}
\author[1]{Daniel W. Crews}
\affil[1]{Computational Plasma Dynamics Lab, University of Washington, Seattle, WA}

\begin{document}
    \maketitle
	
	\begin{abstract}As a method of universal approximation deep neural networks (DNNs) are capable of finding approximate solutions to problems posed
with little more constraints than a suitably-posed mathematical system and an objective function.
Consequently, DNNs have considerably more flexibility in applications than classical numerical methods.
On the other hand they offer an uncontrolled approximation to the sought-after mathematical solution.
This suggests that hybridization of classical numerical methods with DNN-based approximations may be a desirable approach.
In this work a DNN-based approximator inspired by the physics-informed neural networks (PINNs) methodology is used to provide an initial guess
to a Newton's method iteration of a very-high order implicit Runge-Kutta (IRK) integration of a nonlinear system of ODEs, namely the Lorenz system.
In the usual approach many explicit timesteps are needed to provide a guess to the implicit system's nonlinear solver, requiring enough work to make the IRK method infeasible.
The DNN-based approach described in this work enables large implicit time-steps to be taken to any desired degree of accuracy 
for as much effort as it takes to converge the DNN solution to within a few percent accuracy.
This work also develops a general formula for the matrix elements of the IRK method for an arbitrary quadrature order.
\end{abstract}
	
    \section{Introduction}\label{sec:intro}
    Neural networks are universal approximators, so it seems inevitable that they
    take their place in the toolbox of the working applied mathematician and engineer.
    Setting up a neural net approximation for a given problem is fairly easy compared to classical methods of
    discretization.
    The result \textit{looks} like the correct answer to the eye, which for many purposes may be good enough.
    Yet at their current stage of development neural networks seem to provide an \textit{uncontrolled} approximation
    in the sense that there is no obvious relationship between error, iterations, and the nature of the input.
    In contrast, what one might call classical methods, like Newton's root-finding algorithm, are
    a more \textit{controlled} approximation, but require a good guess to get started.
    This project aims to work with the strength of both techniques by hybridizing a neural network with classical numerical analysis.

    Recently, neural networks have been used to solve PDEs in a framework called physics-informed neural networks
    (PINNs)~\cite{pinns}.
    The authors constrained the PINN loss function to satisfy a PDE via automatic differentiation as a continuous
    approximation,
    and experimented with discrete approximations by informing the PINN with
    classical numerical methods like implicit Runge-Kutta (IRK)
    time-steppers.
    However, the output only \textit{looks like} the correct solution, as the neural network loss does not converge
    beyond an L2 error of $\approx 1\%$~\cite{extended_pinn}.
    To solve the problem of convergence, this project uses Newton's method, a beloved root-finding algorithm,
    to correct the neural network prediction for an arbitrarily high order IRK method.
    While Newton's method requires calculation of a system's Jacobian matrix, if one knows physics with which to
    inform a neural network, then one can also calculate the Jacobian of said physics!

    Usually a less-accurate explicit time-step calculates the guess for an implicit solution.
    However, explicit time-steps can be unstable for large steps (\textit{e.g.} the famous CFL condition),
    thereby limiting the implicit step-size.
    Yet the IRK method on Gauss-Legendre points is A-stable, in Butcher's terminology, to any order~\cite{butcher}.
    The implicit time-stepper can take very large steps provided that it has a good guess.
    The advantage of a PINN approach in providing the guess for Newton's method is to overcome the step-size limitation
    of explicitly-predicted stages.

    Section~\ref{sec:lorenz_net} describes the hybrid IRK-integrator method applied to the Lorenz system, which is
    chosen for the chaoticity of its trajectories.
    High-order methods are desirable to ensure that one remains on the correct trajectory when taking a large step.
    Appendix~\ref{sec:irk} develops a closed-form expression for the matrix elements of the IRK method
    to arbitrary order, and Appendix~\ref{sec:newton} then develops the IRK Jacobian for Newton's method.
    If desired, the implementation of all work presented in this report may be found on the author's
    \href{https://github.com/crewsdw/pinns_project}{Github page} using Python and Tensorflow 2.

    \section{Prediction of IRK stages with a neural network and correction}\label{sec:lorenz_net}
    In the following, $m$ denotes the step number, Greek characters $\alpha, \beta, \cdots$ denote vector components, and
    Roman characters $i, j, \cdots$ the RK stages.
    Following Raissi et al., a physics-informed neural net (PINN) can learn the IRK stages given the
    input and an appropriate constraint~\cite{pinns}.
    In the context of the system of ODEs $dy^\alpha/dt=f(\bm{y},t)$ the neural network takes the state $y^{m, \alpha}$ as input
    and for each input has $(n+1)$ outputs, the $n$ RK stages $y^{\alpha}_i$ plus the full step output $y^{m+1, \alpha}$.
    To formulate the loss function, the IRK system for the step $y^{m,\alpha}\to y^{m+1,\alpha}$ is
    rearranged into the following form,
    \begin{align}\label{eq:irk_objective1}
        y^{m,\alpha} &= y^{\alpha}_i - h\sum_{j=0}^{n-1}A^{j}_{i}f_j^\alpha(\bm{y}),\\
        y^{m,\alpha} &= y^{m+1,\alpha} - h\sum_{i=0}^{n-1}\frac{w_i}{2}f_i^{\alpha}(\bm{y})\label{eq:irk_objective2}.
    \end{align}
    The neural net is then trained on its own input $y^{m,\alpha}$ in an autoencoder-like fashion so that the outputs
    $\{y_i^\alpha, y^{m+1,\alpha}\}$ satisfy the system of Eqs.~\ref{eq:irk_objective1} and~\ref{eq:irk_objective2}
    under a mean-square error minimization objective.
    In this way the only data required is that of the initial condition, and not the solution at any later time.
    Requiring only the initial data is a desirable property in an integrator.
    After training to a suitable loss, the network output is given to Newton's method as its initial guess.

    \begin{figure}[h!]
        \centering
    \begin{subfigure}{.45\textwidth}
        \centering
        \includegraphics[width=\linewidth]{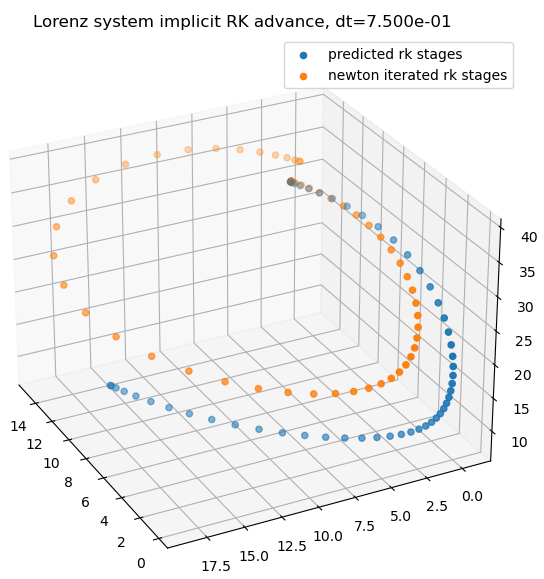}
        \caption{Activation function \texttt{tanh}.}
        \label{fig:sub-first}
    \end{subfigure}
    \begin{subfigure}{.45\textwidth}
        \centering
        \includegraphics[width=\linewidth]{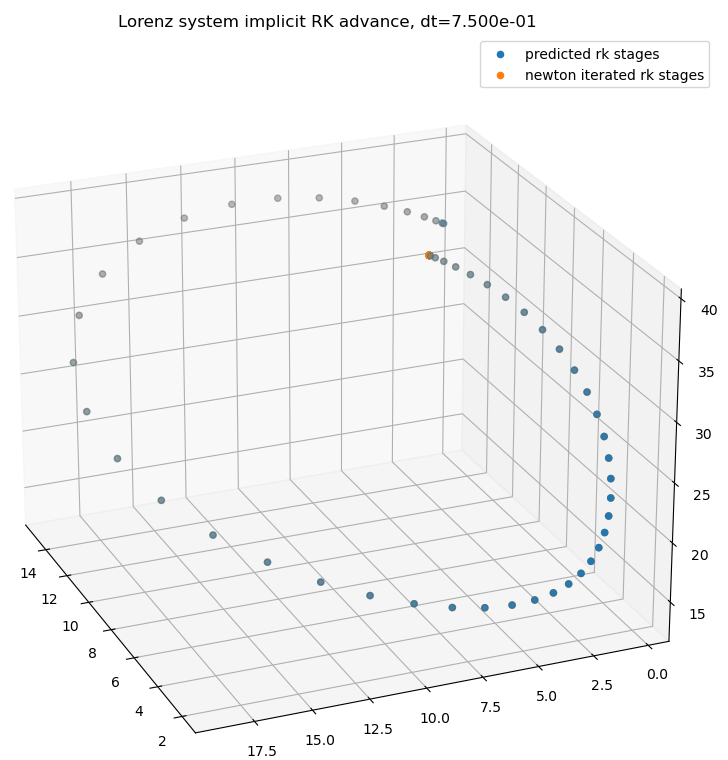}
        \caption{Activation function \texttt{ELU}.}
        \label{fig:sub-second}
    \end{subfigure}
    \caption{An orbit in the Lorenz system integrated using a fifty-stage IRK step of step-size $h=0.75$, solved
    by Newton iteration based on the orbit predicted by a neural network approximation.
    Shown is a comparison of neural net predictions using hyperbolic tangent (\texttt{tanh}) and exponential linear unit
        (\texttt{ELU}) given 10,000 training epochs, and the result of Newton iteration using the neural net output
        as the initial guess. Both guesses allow Newton's method to take over and produce an overall L2 error
        for the implicit system of $10^{-10}$. Clearly, \texttt{ELU} works better here.
        Newton's method is not observed to converge given a worse guess, such as explicit forward Euler to each RK time-stage.}\label{fig:tanh_vs_elu}
    \end{figure}

    \subsection{Application to the Lorenz system}\label{subsec:lorenz_system}
    For an application, consider the famous Lorenz system of ODEs,
    \begin{align}
        \frac{dx}{dt} &= \sigma(y-x),\\
        \frac{dy}{dt} &= x(\rho - z) - y,\\
        \frac{dz}{dt} &= xy - \beta z,
    \end{align}
    with ``original" parameters $\sigma = 10$, $\beta = 8/3$, and $\rho = 28$.
    The Lorenz system's Jacobian matrix is
    \begin{equation}\label{eq:jacobian}
        J^{\alpha}_{\beta} = \begin{bmatrix}-\sigma & \sigma & 0\\ \rho & -1 & -x\\ y & x & -\beta
                                 \end{bmatrix}.
    \end{equation}
    Choosing an initial condition $q_0=[10.54, 4.112, 35.82]$~\cite{wikipedia}, a neural net is trained on
    the initial condition vector of shape $(3, 1)$.
    The net has three hidden layers of three neurons each, and an output layer of shape $(3, n+1)$.
    That such a small number of neurons in each layer ends up working is, to this novice at least, quite remarkable!

    To choose the activation function, Fig.~\ref{fig:tanh_vs_elu} considers fitting the IRK system for a time-step
    of $h = 0.75$ on a 50-stage IRK method using both \texttt{tanh} and \texttt{ELU} activations.
    The function \texttt{ELU} is deemed the winner, at least for this application, and is used subsequently.

    \begin{figure}[h!]
        \centering
    \begin{subfigure}{.45\textwidth}
        \centering
        \includegraphics[width=\linewidth]{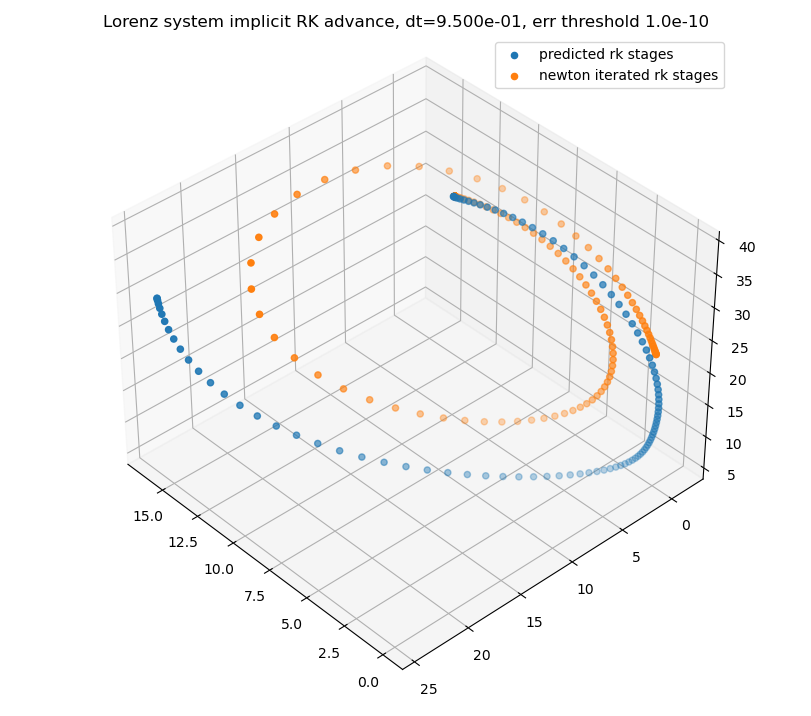}
        \caption{Training epochs: 10,000.}
        \label{fig:sub-first1}
    \end{subfigure}
    \begin{subfigure}{.45\textwidth}
        \centering
        \includegraphics[width=\linewidth]{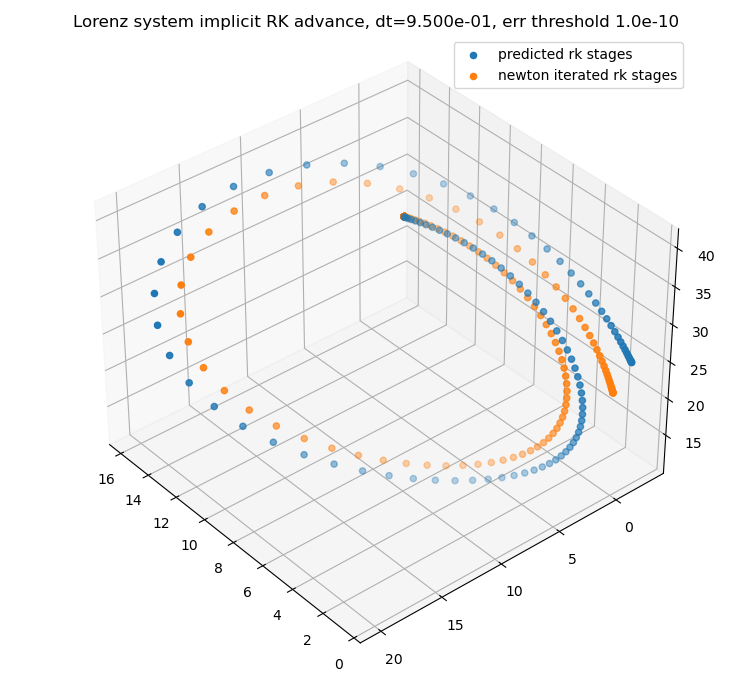}
        \caption{Training epochs: 40,000.}
        \label{fig:sub-second1}
    \end{subfigure}
    \caption{One-hundred stage IRK solution by Newton-iteration of a neural net prediction of 10,000 and 40,000 epochs.
    The Newton solver converges in each case, though a damping factor of $\gamma = 0.25$ is necessary. Convergence
    of the neural net is slow because the optimization algorithm becomes stuck in local minima, possibly that of a
    nearby trajectory. The neural net prediction in each case is good enough for the Newton solver to converge,
        which does not occur if an explicit update is used to guess.
        The Newton solver is much faster than waiting for the neural net to converge, which suggests the utility of
    methods which hybridize ML techniques with classical numerical methods.}\label{fig:0p95_epochs}
    \end{figure}

    \subsubsection{Pushing the limit with 100-stage implicit RK advances}
    Theoretically, a one-hundred stage implicit RK advance on Gauss-Legendre nodes is $\mathcal{O}(h^{200})$-accurate.
    Thus, one can take very large steps and still obtain a quite accurate solution for all times $t\in [t^{m}, t^{m+1}]$
    as defined by the interpolating polynomial of the RK stages.
    However, it is not known by this author to what extent the floating-point precision of the quadrature nodes, weights,
    and RK matrix elements allow machine precision to be actually obtained in the numerical solution.
    This will hopefully be the subject of a future study.

    Increasing the stages to one hundred, and the time-step towards $h=1$, is observed to be difficult for both training
    the neural net and for iterating the net output using Newton's method.
    However, solutions are obtained with sufficient training epochs on the neural net and
    by cranking down the damping factor $\gamma$ on Newton's method.
    A reasonably large time-step is $h=0.80$, with a truncation error of order $\mathcal{O}((0.8)^{200})\approx 10^{-20}$,
    far below machine precision.

    Yet, just to push the limits, Fig.~\ref{fig:0p95_epochs} investigates the case $h=0.95$.
    The neural net did converge to the correct solution, but required a large number of steps.
    Even with the good guess provided by the neural net,
    a significant damping factor of $\gamma=0.25$ was necessary for Newton's method to converge at this large step.
    A 100-stage step of step-size $h=1$ was also obtained as an experiment, but required a large amount of computational effort.
    Further, any $h>1$ should amplify the error rather than decrease it, so this was not investigated.

\begin{figure}
        \centering
        \includegraphics[width=0.5\textwidth]{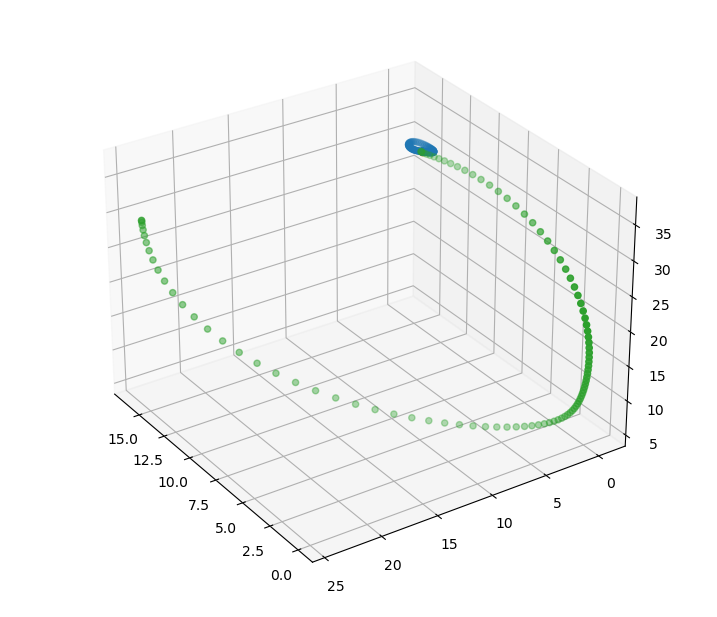}
        \caption{The neural net objective is to linearly combine the output stages through the RK matrix to match
        the initial condition. This figure visualizes the technique for a 100-stage IRK method which is stuck in a local
            minimum, with the output stages
            shown in green and the objective points in blue.
        The method is convergent when the blue points all match the initial condition.}\label{fig:net_out}
    \end{figure}

    Figure~\ref{fig:net_out} shows the uncorrected neural net outputs along with their combination
    into the objective function (the initial condition) as an aid to visualize the method.
    For illustration, this is shown for a case which is slowly converging, as it is stuck in a local minimum.
    When converged, all the blue dots will coincide with the initial condition $q_0$.

\begin{figure}[h!]
        \centering
        \begin{subfigure}{.45\textwidth}
            \centering
            \includegraphics[width=\linewidth]{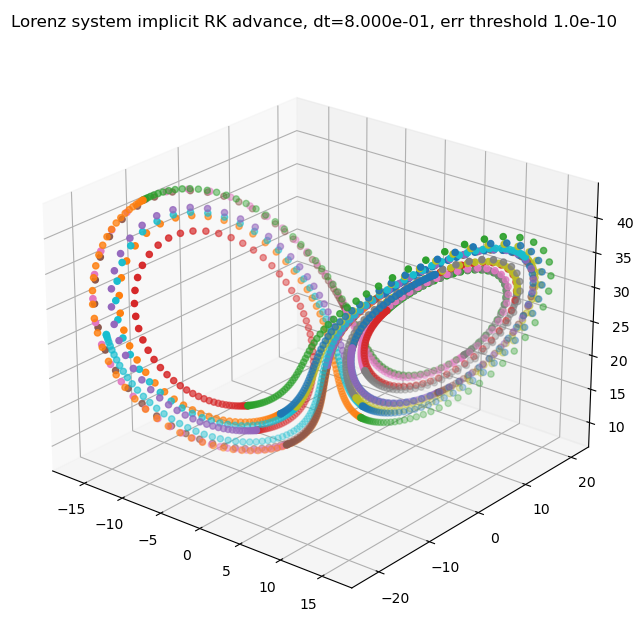}
            \caption{Perspective on Lorenz attractor.}
            \label{fig:sub-first2}
        \end{subfigure}
        \begin{subfigure}{.45\textwidth}
            \centering
            \includegraphics[width=\linewidth]{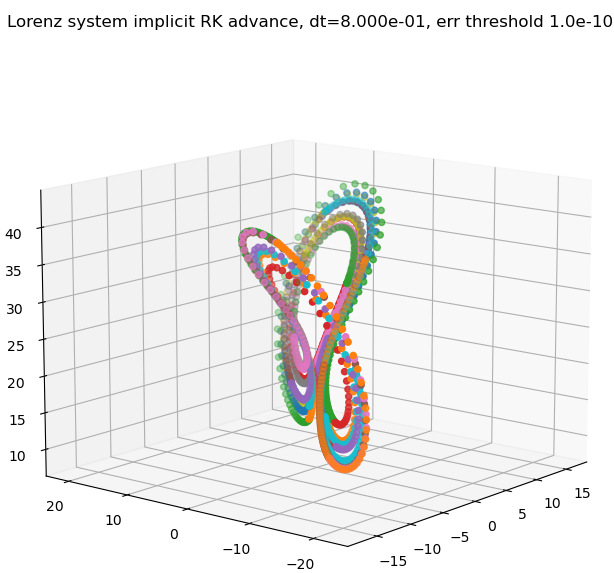}
            \caption{Another attractor perspective.}
            \label{fig:sub-second2}
        \end{subfigure}
        \caption{A trajectory of the Lorenz system integrated with ten time-steps of a 100-stage implicit RK Gauss-Legendre
        method of step-size $h=0.8$, solved by Newton iteration from a neural net prediction of the RK stages.
        Each color shown corresponds to the stages of a single step.}\label{fig:lorenz_attractor}
    \end{figure}

    Finally, the neural net predictor-Newton's method corrector 100-stage IRK time-stepper was applied to multiple
    time-steps.
    Figure~\ref{fig:lorenz_attractor} shows the Lorenz attractor traced by the trajectory integrated out to ten steps
    each of step-size $h=0.8$.
    While the IRK method with step-size $h=0.8$ has theoretical truncation error $\mathcal{O}(10^{-20})$,
    the Newton iteration L2 error threshold is set to $10^{-10}$. 

    \section{Conclusions and extensions}\label{sec:conclusions}
    Neural nets have immense promise as a universal approximation system,
    broadened by their being relatively easy to set-up for a wide variety of problems
    with automatic differentiation baked-in.  
    However, the objective function landscape is full of local minima.
    While the result of training a neural net generally \textit{looks like} the correct answer,
    there are usually convergence issues due to lack of constraints in the posed optimization problem.

    This project combined the neural net output with an approximation of classical numerical analysis,
    namely Newton's method.
    The neural network prediction is used as a guess in order to solve the nonlinear system associated with an arbitrarily
    high-order implicit time-integration of a system of ODEs.
    Assisted by the good guesses of the neural network, Newton's method was found to converge for an
    $\mathcal{O}(h^{200})$-accurate method.

    While this study looked at the Lorenz system, there's no reason that this ODE couldn't correspond to the semi-discrete
    ODE system of a spatially-discretized PDE\@.
    The only difficulty of such extensions is to calculate the system Jacobian, which is shown in Appendix~\ref{sec:newton}
    to consist of the Runge-Kutta coefficient matrix and the Jacobian of the ODE system.
    As Appendix~\ref{sec:irk} establishes a closed-form expression for the RK coefficient matrix, to apply the method to
    any system of ODEs one need only calculate its Jacobian.
    Further extensions of this project can look at more sophisticated iterative methods which incorporate information
    on second derivatives, such as BFGS.
	
	\section{Acknowledgements}
	The author would like to thank J. Bakarji, S.L. Brunton, J.N. Kutz, U. Shumlak, and A.D. Stepanov for teaching and helpful discussions.

    \bibliography{main}
    \bibliographystyle{unsrt}

    \newpage
    \appendixpage
    \appendix

        \section{Implicit RK integration using Gauss-Legendre quadrature}\label{sec:irk}
    This section considers the implicit Runge-Kutta (RK) method using Gauss-Legendre (GL) quadrature and develops
    a closed-form expression for the RK coefficient matrix.
    First, recall that
    \begin{equation}\label{eq:gl_quad}
        \int_{-1}^{1}f(\xi)d\xi \approx \sum_{i=0}^{n-1}w_{i}f(\xi_i)
    \end{equation}
    with $\{w_i\}_{i=0}^{n-1}$ quadrature weights and $\{\xi_i\}_{i=0}^{n-1}$ the roots of $P_n(\xi)$,
    the $n$'th Legendre polynomial.
    The integration is exact for polynomials of order $\leq 2n-1$.
    There is a beautiful theory of the interpolation polynomials through GL quadrature points;
    it is best understood in terms of the Legendre polynomial completeness relation~\cite{hassani}
    \begin{equation}\label{eq:spherical_completeness}
        \sum_{s=0}^{\infty}\Big(s + \frac{1}{2}\Big)P_s(x)P_s(y) = \delta(y-x),\quad\quad -1\leq x,y\leq 1.
    \end{equation}

    \begin{thm}
        The interpolation polynomials $\ell_j(\xi)$ through $n$ Gauss-Legendre nodes satisfy
    \begin{equation}\label{eq:lagrange_exp}
        \ell_j(\xi) = w_j\sum_{s=0}^{n-1}\Big(s + \frac{1}{2}\Big)P_s(\xi_j)P_s(\xi),\quad\text{ and }\quad \ell_j(\xi_i) = \delta_{ij}
    \end{equation}
    with $\{w_j\}_{j=0}^{n-1}$, $\{\xi_j\}_{j=0}^{n-1}$ the quadrature weights and nodes respectively\footnote{Considering the
    summation $\lim_{j\to\infty}\sum_{j=0}^{n-1} \ell_j(\xi)f(\xi)$ suggests that the
    the interpolation polynomial plays a role in discrete integration as if it were a $\delta$-sequence on $[-1,1]$ in the
    sense of distributions~\cite{stakgold}.}.\\
        \vspace{0.25cm}

        \textbf{Proof}:
        The interpolation polynomials are defined in Lagrange form as
        \begin{equation}\label{eq:lagrange}
            \ell_j(\xi) = \prod_{\substack{k=0\\ k\neq j}}^{n-1}\frac{\xi - \xi_k}{\xi_j-\xi_k}.
        \end{equation}
        Now, note that in addition to continuous orthogonality
        \begin{equation}\label{eq:cont_orth}
            \int_{-1}^1 P_k(\xi)P_{\ell}(\xi)d\xi = \frac{1}{k+\frac{1}{2}}\delta_{k\ell}
        \end{equation}
        the Legendre polynomials are discretely orthogonal~\cite{dunkl}~\cite{bulirsch},
        \begin{equation}
            \sum_{s=0}^{n-1}w_s P_k(\xi_s)P_{\ell}(\xi_s) = \frac{1}{k+\frac{1}{2}}\delta_{k\ell}\label{eq:gauss_orthogonal}
        \end{equation}
        with $w_s$ the quadrature weights.
        This follows from the order $2n-1$ quadrature property as the sum corresponds to the continuous integral $\int_{-1}^{1}P_k(x)P_\ell(x)dx$.
        Expansions in the Lagrange basis of Eqn.~\ref{eq:lagrange} and in the Legendre basis are linearly related
        due to the interpolation property $\ell_j(\xi_i) = \delta_{ij}$,
        \begin{equation}
            f(\xi) = \sum_{j=0}^{n-1}f_j\ell_j(\xi) = \sum_{k=0}^{n-1}c_{k}P_k(\xi),\quad\implies\quad f_j = \mathcal{V}_j^{k}c_k\label{eq:vandermonde1}
        \end{equation}
        where $\mathcal{V}^k_j = P_k(\xi_j)$ is termed the (generalized) Vandermonde matrix.
        The inverse transform follows directly from Eqn.~\ref{eq:gauss_orthogonal} as
        $(\mathcal{V}^{-1})_j^k=w_j P_k(\xi_j)(k+\frac{1}{2})$.
        In particular, the Lagrange interpolant itself is
        \begin{equation}
            \ell_j(\xi) = w_j\sum_{s=0}^{n-1}\Big(s + \frac{1}{2}\Big)P_s(\xi_j)P_s(\xi) \approx w_j\delta(\xi-\xi_j).\label{eq:lagrange_exp2}
        \end{equation}
    \end{thm}
    \subsection{The Runge-Kutta method for Gauss-Legendre quadrature}\label{subsec:rk_method}
    The IRK coefficient matrix for GL stage-points
    is developed to arbitrary order using the interpolation polynomial of Eq.~\ref{eq:lagrange_exp}.
    To recall the Runge-Kutta method, consider the initial value problem
    \begin{equation}\label{eq:rk}
        \frac{dy}{dt} = f(t,y),\quad y(t_m) \equiv y^m.
    \end{equation}
    The journey to a new step $y^m\to y^{m+1}$ with step-size $h$ is calculated via a sequence of $n$ stages
    $\{y^{m, c_i}\}_{i=0}^{n-1}$ at the stage times $\{t_i\}_{i=0}^{n-1}$ according to the system
    \begin{align}
        y^{m, c_i} &= y^m + h\sum_{j=0}^{n-1} A_i^{j}f(y^{m, c_j}),\\
        y^{m+1} &= y^m + h\sum_{i=0}^{n-1} b_{i}f(y^{m, c_i})
    \end{align}
    where $A_i^j$ is the RK coefficient matrix and $\{b_i\}_{i=0}^{n-1}$ the stage combination coefficients.
    The stage locations $\{t_i\}$, combination coefficients $\{b_i\}$, and RK matrix $A_i^j$ are said to form the Butcher tableau~\cite{butcher}.
    \begin{thm}
        The Butcher tableau for the Gauss-Legendre IRK method consists of:
        \begin{itemize}
            \item Stage times: the GL nodes $\{\xi_i\}_{i=0}^{n-1}$ transformed from $[-1, 1]\to [t^m, t^{m+1}]$ according to
            \begin{equation}\label{eq:affine}
                \xi = \frac{h}{2}(t - \bar{t}^m)
            \end{equation}
            with step-size $h \equiv t^{m+1} - t^m$ and mid-point $\bar{t}^m \equiv \frac{1}{2}(t^{m+1} + t^m)$,
            \item Combination coefficients: the GL quadrature weig                                                                                                                  hts, $b_i \equiv w_i/2$ (with $\frac{1}{2}$ from the transform).
            \item RK coefficient matrix: the matrix elements are given by
            \begin{equation}\label{eq:rk_matrix}
                A_i^j = \frac{w_j}{2}\Big(1 + \sum_{s=0}^{n-1}P_s(\xi_j)\frac{P_{s+1}(\xi_i) - P_{s-1}(\xi_i)}{2}\Big)
            \end{equation}
            where $P_{-1}(\xi)\equiv 1$, and by definition of the quadrature nodes $P_{n}(\xi_j)=0$.
        \end{itemize}

        \textbf{Proof}:
        Considering the ODE $\frac{dy}{dt} = f(t,y)$, project the RHS onto the basis defined by the interpolants $\ell_i(t)$
        through the GL quadrature points $\{t_i\}_{i=1}^n$ defined by the affine transform of Eq.~\ref{eq:affine},
        \begin{equation}\label{eq:interpolated_rhs}
            f(y, t) \approx \sum_{j=0}^{n-1} f(y^{m, t_j})\ell_j(t).
        \end{equation}
        Solving the ODE by integration, one has
        \begin{align}
            y(t^{m, c_i}) - y(t^m) &= \int_{t^m}^{t^{m, c_i}} f(y,t)dt\\
                                &\approx \sum_{j=0}^{n-1}f(y^{m, t_j})\int_{t_m}^{t_{m, c_i}}\ell_j(t)dt
        \end{align}
        so the calculation reduces to integration of the interpolating polynomial.
        Using Eq.~\ref{eq:lagrange_exp} and the affine transform, this further reduces to integrating the Legendre polynomial,
        \begin{equation}\label{eq:int_legendre}
            \int_{t^m}^{t^{m, c_i}}\ell_j(t)dt = \frac{h}{2}w_{j}\sum_{s=0}^{n-1}\Big(s+\frac{1}{2}\Big)P_s(\xi_j)\int_{-1}^{\xi_i}P_s(\xi)d\xi.
        \end{equation}
        Recalling the identity
        \begin{equation}\label{eq:legendre_identity}
            P_s(\xi) = \frac{1}{2s+1}\frac{d}{d\xi}\Big(P_{s+1}(\xi) - P_{s-1}(\xi)\Big)
        \end{equation}
        the Legendre polynomial integrates to
        \begin{equation}\label{eq:int_legendre2}
            \int_{-1}^{\xi_i}P_s(\xi_i)d\xi = \frac{1}{2}\frac{1}{s+\frac{1}{2}}\Big[(P_{s+1}(\xi_i) - P_{s-1}(\xi_i)) -
        (P_{s+1}(-1) - P_{s-1}(-1))\Big].
        \end{equation}
        The boundary term $\xi = -1$ cancels for all $s>0$, but for $s=0$ an extra constant factor is picked up.
        This establishes Eq.~\ref{eq:rk_matrix}.
        Lastly, the combination coefficients are found by integrating all the way to $t^{m+1}$, for an integral
        $\int_{-1}^1 P_s(\xi)d\xi = 2\delta_{s0}$ by orthogonality.
        Thus,
        \begin{equation}\label{eq:stage_combination}
            y^{m+1} = y^m + h\sum_{i=1}^n \frac{w_{i}}{2}f(y^{m, c_i}).
        \end{equation}
    \end{thm}
    The class of RK methods is very wide, and a large body of research exists on developing general methods.
    One might get the impression that high-order methods are difficult to construct.
    However, the preceding shows that when Gauss-Legendre points are utilized the method may be developed
    to arbitrarily high order, where only the quadrature nodes and weights are needed.
    For this project methods up to $n=100$ are acquired using an online calculator~\cite{gauss:calculator}, though
    recent research has developed explicit formulas via asymptotic methods to calculate the nodes and weights to
    floating point precision for any $n\geq 20$~\cite{townsend}~\cite{bogaert}.

    \section{The Newton-Raphson method for vector systems}\label{sec:newton}
    The Newton-Raphson root-finding method is a beloved approximation and is close to the hearts of many.
    To quickly review, given a system of equations in vector form
    \begin{equation}\label{eq:vector_equation}
        f^{\alpha}(\bm{x}) = 0
    \end{equation}
    its linearization about a point $x_0^\alpha$ is, using the Jacobian matrix $J^\alpha_{\beta} \equiv \frac{\partial f^\alpha}{\partial x^\beta}$,
    \begin{equation}\label{eq:linearization}
        f^{\alpha}(\bm{x}_0) + J^\alpha_\beta(\bm{x}_0)(x^\beta - x_0^\beta) = 0
    \end{equation}
    where summation convention is used.
    Based on an initial guess $x_0^\beta$, one solves the linear system,
    \begin{equation}\label{eq:newton_raphson}
        J^\alpha_\beta(\bm{x}_0)(x^\beta - x_0^\beta) = -f^{\alpha}(\bm{x}_0)
    \end{equation}
    for unknown $q^\beta \equiv x^\beta - x_0^\beta$.
    An iteration is then obtained as
    \begin{equation}\label{eq:iterate}
        x_1^\beta =  q^\beta + x_0^\beta.
    \end{equation}
    One then repeats the process, sometimes adding a damping factor $|\gamma|\leq 1$ to the iteration to ensure convergence.
    Evidently, to use the method one must find the Jacobian matrix of the system.

    \subsection{Newton-Raphson iteration for the implicit Runge-Kutta method}\label{subsec:newton_irk}
    Considering a vector unknown $y^\alpha(t)$, the RK method integrates the system of ODEs
    \begin{equation}\label{eq:vector_ode}
        \frac{dy^\alpha}{dt} = f^\alpha(\bm{y}, t).
    \end{equation}
    Denoting the right-hand side evaluation for the RK stages as the vector $\bm{k}^i \equiv k^{\alpha, i}$,
    for a vector unknown the RK system may be written as
    \begin{align}
        y^{m+1, \alpha} &= y^{m, \alpha} + h b_{i}k^{\alpha, i},\\
        k^{i, \alpha} &= f^{\alpha}(t^m + c_ih, y^m + h A_i^j k^\alpha_j)
    \end{align}
    with RK matrix $A_i^j$ and coefficients $c_i$, $b_i$ found in Appendix~\ref{subsec:rk_method}.
    The second of these equations forms the implicit system of equations to be solved, so it's sufficient to apply
    Newton's method to
    \begin{equation}\label{eq:equals_zero}
        k^{i, \alpha} - f^{\alpha}(\bm{k}^i) = 0.
    \end{equation}
    Having linearized, for each iteration one solves the system
    \begin{equation}\label{eq:irk_newton}
        \mathbb{J}^{i,\alpha}_{j,\beta} q^{j, \beta} = f^{\alpha}(\bm{k}^i) - k^{i, \alpha}.
    \end{equation}
    The iteration update is then $k \texttt{ += } \gamma q$ with $\gamma$ an iterate damping factor.
    The Jacobian $\mathbb{J}$ is a linear operator
    $\mathbb{J}:\mathbb{R}^{\Lambda, n}\to \mathbb{R}^{\Lambda, n}$ for $\Lambda$ vector components and the $n$ RK stages.
    As a tensor-like object it is found most easily in index notation.
    \begin{thm}
        The Jacobian tensor $\mathbb{J}^{i, \alpha}_{j, \beta}$ of the implicit Runge-Kutta method is
        \begin{equation}\label{eq:irk_jacobian}
            \mathbb{J}^{i,\alpha}_{j,\beta} = \delta^i_j\delta^\alpha_{\beta} - h A^i_j J^{\alpha, i}_{\beta}
        \end{equation}
        where $J^{\alpha, i}_{\beta} \equiv \partial_{y^\beta}f^\alpha(\bm{k}^i)$ is the Jacobian matrix of the
        original ODE evaluated at the i'th RK stage.\\

        \textbf{Proof}: Using the symbol $\partial_{k^{j, \beta}}$, one has
        \begin{align}
            \mathbb{J}^{i,\alpha}_{j,\beta} &= \partial_{k^{j, \beta}}(k^{i, \alpha} - f^{\alpha}(\bm{k}^i))\\
                                            &= \delta^i_j\delta^\alpha_{\beta}
                                                                            -
            \partial_{k^{j, \beta}}f^{\alpha}\Big(\bm{k}(\bm{y})\Big)\\
                    &= \delta^i_j\delta^\alpha_{\beta} - (\partial_{k^{j, \beta}}y^{i, \gamma})\partial_{y^{i,\gamma}}f^{\alpha}(\bm{y})\\
                    &= \delta^i_j\delta^\alpha_{\beta} - h A^i_j \delta_{\beta}^{\gamma} J^{\alpha, i}_{\gamma}\\
                    &= \delta^i_j\delta^\alpha_{\beta} - h A^i_j J^{\alpha, i}_{\beta}
        \end{align}
    \end{thm}
    In summary, the Jacobian for an IRK method of arbitrarily high order can be calculated using: i) the RK coefficient
    matrix, and ii) the Jacobian of the original ODE system.
    The products may be done numerically using \texttt{numpy.einsum()}, for example.
    Of course, one needs a good guess for the RK stages $y^{\alpha, i}$ to get started.
    A guess for large time-steps can be furnished with a neural network.

\end{document}